\documentclass[12pt]{article}

\usepackage{fullpage}
\usepackage{amsmath}
\usepackage[pdftex]{graphicx}
\usepackage{graphicx}
\usepackage{amssymb}
\usepackage{epstopdf}
\usepackage{booktabs}
\usepackage{hyperref}
\usepackage{cite}
\usepackage[font={small}, margin=1cm]{caption}
\usepackage[titletoc,title]{appendix}

\newtheorem{prop}{Proposition}
 \title{Three-Dimensional Multi-Tethered Satellite Formation with the Elements Moving Along Lissajous Curves }

\author{D. Yarotsky\footnote{Institute for Information Transmission Problems, Russian Academy of Sciences, Bolshoy Karetny per. 19, 127051 Moscow, Russian Federation}  \and  
V. Sidorenko\footnote{Keldysh Institute of Applied Mathematics, Russian Academy of Sciences, Miusskaya Sq., 4, 125047 Moscow, Russian Federation} \and
D. Pritykin\footnote{Moscow Institute of Physics and Technology, Institutskiy per., 9, 141700, Dolgoprudny, Moscow Region, Russian Federation}}
\date{}				

\begin{document}
\bibliographystyle{plainnat}
 \maketitle
 \begin{abstract}

This note presents a novel approach to maintain three-dimensional multi-tethered satellite formation in space. For a formation consisting of a main body connected by tethers with several deputy satellites (the so-called ``hub-and-spoke'' configuration) we demonstrate that under proper choice of the system's parameters the deputy satellites can move along Lissajous curves in the plane normal to the local vertical with all tethers stretched; the total force due to the tension forces acting on the main satellite is balanced in a way allowing it to be in relative equilibrium strictly below or strictly above the system's center of mass. We analyze relations between the system's essential parameters and obtain conditions under which the proposed motion does take place. We also study analytically the motion stability for different configurations and whether the deputy satellites can collide or the tethers can entangle. Our theoretical findings are corroborated and validated by numerical experiments.


\end{abstract}

\section{Introduction}\label{intro}

Three-dimensional satellite formations are often discussed in connection with multi-point measurements needed for atmospheric, geodetic or plasma physics studies. To simplify control strategies and to minimize fuel consumption, tethers can be used to maintain desired relative positions of satellites in the formation flying. For the first time three-dimensional multi-tethered formations were discussed probably by \cite{B1983}, who proposed double-pyramid configurations. It seems that the most straightforward way to obtain a multi-tethered formation is to deploy from the main satellite several tethers with deputy satellites at their ends. To specify such formations \cite{P2008} introduced the term ``hub-and-spoke''. Behavior of ``hub-and-spoke'' multi-tethered formations has been studied for different dynamical environments:  in circular orbit \cite{A2013,A2015}, in elliptic orbit \cite{A2014}, in halo-orbit \cite{Z2008, C2014} and near collinear Lagrangian points \cite{W2008}.

To keep the tethers taut the combination of rotation with gravity-gradient forces is usually proposed. Among other opportunities the relatively new concept of the Tethered Coulomb Structure (TCS) is worth mentioning \cite{S2009, S2011, P2012}. In this case the satellites are electrostatically charged to produce repulsive forces between them. Nevertheless, it looks as if Coulomb repulsive forces can be used to prevent the slack of short enough tether: in \cite{S2009, S2011, P2012} the discussed length is 10~m by order of magnitude.

\begin{figure}[htb]
\begin{center}
\includegraphics[width=0.5\textwidth,keepaspectratio]{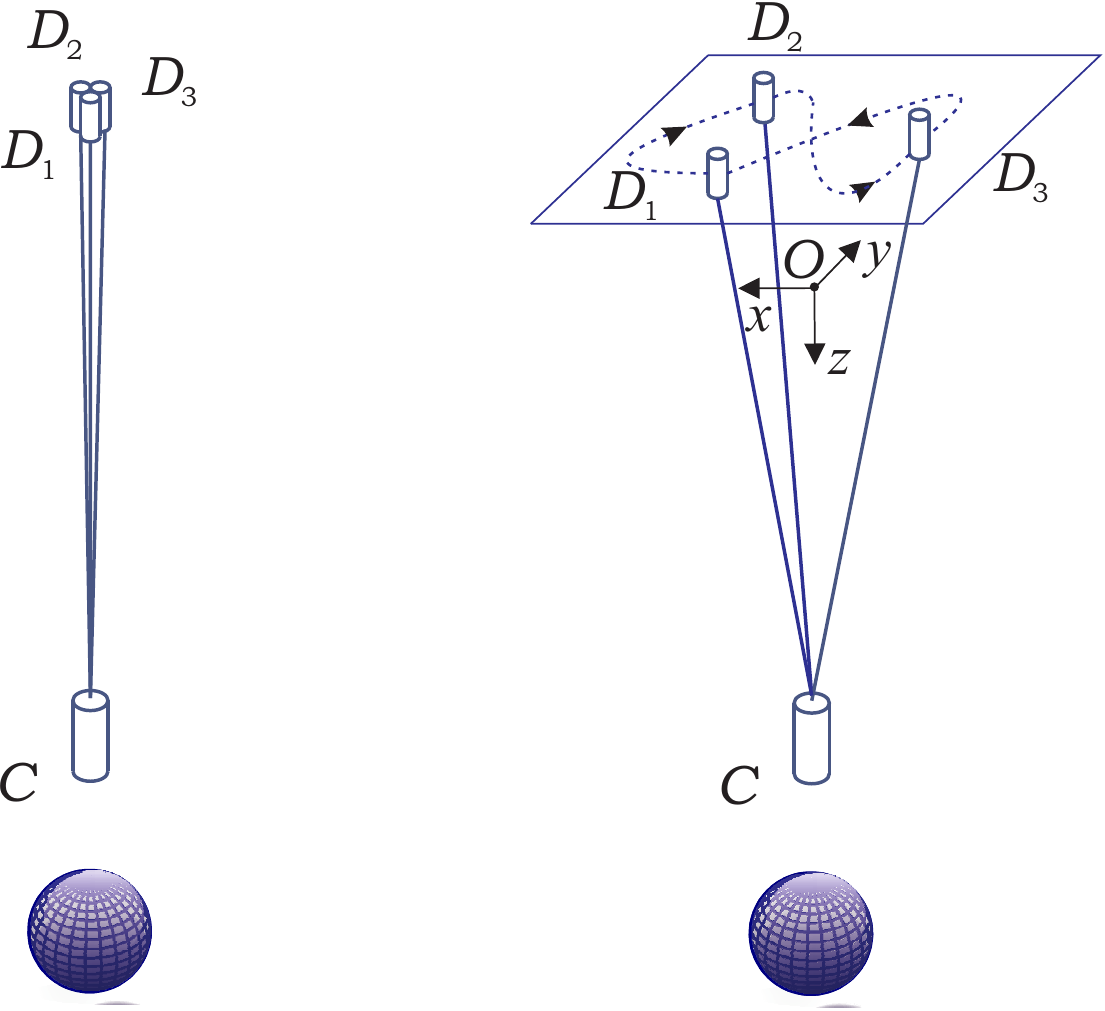}
\end{center}
    \caption{Multi-tethered formation consisting of the main body $C$ and three deputy satellites $D_1, D_2, D_3$. On the left the system is in relative equilibrium with all tethers aligned along the local vertical. On the right the deputy satellites move along a Lissajous curve (see also the animation provided as
Electronic Supplementary Material).}
    \label{fig:fig1}
\end{figure}

To give an idea of our approach we begin with the system of two bodies connected by a single tether aligned along the local vertical; the mass center of the system moves in a circular orbit. Small oscillations of this system around local vertical are a combination of in-plane and out-of-plane natural oscillations with incommensurable frequencies \cite{C2008,S2010}. The motion of end bodies in these oscillations can be roughly described as a motion along curves densely filling certain areas on planes normal to the local vertical. Then let us consider the degenerate ``hub-and-spoke'' configuration in the relative equilibrium with all tethers aligned along the vertical (Fig.~\ref{fig:fig1}, left) The displacement of a single deputy satellite from the relative equilibrium position causes oscillations whose frequencies differ from those inherent in the preceding case. With the proper tuning of the system's parameters these frequencies can be made commensurable resulting in the motion of the deputy satellite along a Lissajous curve in the plane, normal to the local vertical (Fig.~\ref{fig:fig1}, right).  Naturally all other deputy satellites can also be put in motion along similar curves. The choice of initial conditions allows to avoid collisions among them and to ensure the balance of the tension forces applied to the main satellite so as to preserve its relative equilibrium.

We suppose that the described structure can be useful for some applications or at least become the starting point for the development of new approaches to maintain 3D multi-tethered satellite formations in space.

In Section 2 we start our study with a simplified dynamical model of multi-tethered formation (point masses + weightless tethers). In Section 3 we consider the linearized dynamics near vertical equilibrium and describe the oscillations in the system. In Section 4 we consider deputy satellite formations moving along Lissajous curves. In Section 5 we present the results obtained by numerical simulation of the system's dynamics.

\section{Deputy satellite dynamics at small deviations from the relative equilibrium }
As mentioned in the Introduction, we will consider the multi-tethered satellite formation comprising $N + 1$ bodies: the main satellite $C$ of mass $m_C$ and $N$ deputy satellites $D_1, ..., D_N$ (each of mass $m_D$) linked to the main satellite (but not to each other) by identical extensible tethers; the tethers' masses are ignored.

To write down the equations of motion, we shall introduce a Local Vertical Local Horizontal (LVLH) reference frame $Oxyz$, centered on the position of the system's center of mass (CoM) in its nominal orbital motion: $Oz$ axis is aligned with the local vertical and oriented towards Earth's center, $Ox$ runs tangentially to the orbit in the direction of the CoM motion, and $Oy$ axis is directed along the normal line to the orbit plane (Fig.~\ref{fig:fig1}). It is assumed that nominally the system's center of mass moves in circular orbit with the mean motion $\omega_0$. The unit vectors in the directions of the axes $x,y,z$ will be denoted as ${\mathbf e}_x,
{\mathbf e}_y, {\mathbf e}_z$ respectively.

In the LVLH frame the motion of any satellite in the considered formation can be described by the Hill-Clohessy-Wiltshire (HCW) equations \cite{C1960}:

\begin{align}
\ddot{x} &= 2\omega_0\dot{z}+T_x/m,\notag\\
\ddot{y} &= \omega_0^2 y+T_y/m,\\
\ddot{z} &= -2\omega_0\dot{x}+3\omega_0^2z+T_z/m,\notag
\end{align}
where $T_x, T_y, T_z$ denote the components of the (total) tether tension $\mathbf{T}$ applied to a given body of mass $m$.

We shall adopt the usual visco-elastic model of massless tether, hence the tension force applied to the $i$th deputy satellite is
\begin{equation}
\mathbf{T}_i=\mathbf{1}_{(|\mathbf{r}_C-\mathbf{r}_i |>l_0 )}\cdot\left[k(|\mathbf{r}_C-\mathbf{r}_i |-l_0 )+b \frac{d}{dt} |\mathbf{r}_C-\mathbf{r}_i |\right]\frac{(\mathbf{r}_C-\mathbf{r}_i)}{|\mathbf{r}_C-\mathbf{r}_i|},
\end{equation}
where $l_0$ is the slack tether length, $k$ the elastic coefficient, $b$ the damping coefficient, and
\[
    \mathbf{1}_{(|\mathbf{r}_C-\mathbf{r}_i |>l_0 )}=
\begin{cases}
    1, & \text{if } |\mathbf{r}_C-\mathbf{r}_i |>l_0\\
    0,              & \text{otherwise}
\end{cases}
\]

The vertical equilibrium of the system takes place if all tethers rest ($\mathbf{\dot{r}}_C=\mathbf{\dot{r}}_i=0$) aligned along the local vertical:
\begin{align}\label{eq:equi}
x_i^* &=x_C^*, \quad y_i^* =y_C^*=0, \\
z_i^* &=-\frac{m_C}{N m_D}z_C^*={l_0}\left(\frac{N m_D+m_C}{m_C} -3\frac{m_D\omega_0^2}{k}\right)^{-1}. \notag
\end{align}

Here and below equilibrium quantities are marked by the asterisk *. In contrast to the $y$ and $z$ coordinates, the $x$ coordinate in this configuration, though shared by all satellites, is arbitrary, since the equations are invariant with respect to translations along the orbit (i.e., along the $x$ axis in the approximation of the orbital dynamics, provided by the HCW equations).

For the above configuration to be valid the denominator in the last formula must be positive, i.e. the tethers must be sufficiently rigid to counteract the microgravity:

\begin{equation}\label{eq:stab1}
k>3{\omega }^2_0\frac{m_Cm_D}{Nm_D+m_C}.
\end{equation}

Denoting the tether's length in the equilibrium configuration by $l_*=|{{\mathbf r}}^{*}_{C} - {{\mathbf r}}^*_i|$, we obtain a useful relation

\begin{equation} \label{eq:len}
\frac{l_* - l_0}{l_*}=\frac{3{\omega }^2_0}{k}\frac{m_Cm_D}{Nm_D+m_C}.
\end{equation}

Denote
\[\lambda_* = \frac{l_* - l_0}{l_*}=\frac{3{\omega }^2_0}{k}m_r, \quad m_r = \frac{m_Cm_D}{Nm_D+m_C}.\]

The tension forces in all the tethers in the equilibrium configuration (\ref{eq:equi}) have the same value ${\mathbf T}^*$ given
by the obvious formula
\begin{equation}
{\mathbf T}^*= k(l_* - l_0){\mathbf e}_z.
\end{equation}

It is not difficult to derive an approximate expression for the tension forces in case of
small displacements of the main satellite and the $i$th deputy satellite with respect to equilibrium (\ref{eq:equi}):
\begin{align*}
{\mathbf T}_i = & {\mathbf T}^* + k\lambda_*(\Delta x_C - \Delta x_i){\mathbf e}_x + k\lambda_*(\Delta y_C - \Delta y_i){\mathbf e}_y + \\
&+\left[k(\Delta z_C - \Delta z_i) + b(\Delta \dot{z}_C - \Delta \dot{z}_i)\right]{\mathbf e}_z.
\end{align*}

Here $\Delta x_C$, $\Delta y_C$, $\Delta z_C$ are the components of the main satellite's displacements, and
$\Delta x_i$, $\Delta y_i$, $\Delta z_i$ are the components of the $i$th deputy satellite displacements.

With the linearized tension the HCW equations read
\begin{align}
\Delta\ddot{x_i} &= 2\omega_0 \Delta\dot{z_i} + \lambda_*\frac{k}{m_D}(\Delta x_C -\Delta x_i),\notag \\
\Delta\ddot{y_i} &= -\omega_0^2 \Delta y_i + \lambda_*\frac{k}{m_D}(\Delta y_C - \Delta y_i), \\
\Delta\ddot{z_i} &= -2\omega_0 \Delta\dot{x_i} +3\omega_0^2 \Delta z_i + \frac{k}{m_D}(\Delta z_C-\Delta z_i) +\frac{b}{m_D} (\Delta \dot{z}_C - \Delta \dot{z}_i) \notag
\end{align}
for the deputy satellites, and
\begin{align}
\Delta\ddot{x_C} &= 2\omega_0 \Delta\dot{z_C} - \lambda_*\frac{k}{m_C}\sum^N_{i=1}(\Delta x_C -\Delta x_i),\notag \\
\Delta\ddot{y_C} &= -\omega_0^2 \Delta y_C - \lambda_*\frac{k}{m_C}\sum^N_{i=1}(\Delta y_C -\Delta y_i), \\
\Delta\ddot{z_C} &= -2\omega_0 \Delta\dot{x_C} +3\omega_0^2 \Delta z_C - \frac{k}{m_C}\sum^N_{i=1}(\Delta z_C-\Delta z_i) - \frac{b}{m_C}\sum^N_{i=1}(\Delta \dot{z}_C - \Delta \dot{z}_i) \notag
\end{align}
for the main satellite.

\section{Decoupling the equations of motion}
The system of equations derived above can be split up into three independent groups by taking appropriate linear combinations, separately for $x$,  $y$ and $z$ components:

\begin{enumerate}
\item  By taking the sum of the equation for the main satellite with weight $m_C$ and all the respective equations for deputy satellites with weights  $m_D$ we obtain a triple of scalar equations for
\[\frac{1}{Nm_D+m_C}\left(m_C\Delta {{\mathbf r}}_C + m_D\sum^N_{i=1}{\Delta {{\mathbf r}}_i}\right),\]
i.e., for the motion of the whole system's center of mass.

\item  By taking the sum of the equations for deputy satellites with coefficients $1/N$ and subtracting from it the equation for the main satellite we obtain a triple of scalar equations for
\[\frac{1}{N}\sum^N_{i=1}{{\Delta {\mathbf r}}_i}-{\Delta {\mathbf r}}_C,\]
i.e. describing the relative motion between the main satellite and the auxiliary satellites' center of mass.

\item  Finally, if we take any linear combination of equations for deputy satellites with some coefficients $\sigma_i$ such that $\sum^N_{i=1}{\sigma_i}=0$, we obtain a triple of scalar equations for $\sum^N_{i=1}{\sigma_i}\Delta {{\mathbf r}}_i$, which can be interpreted as a partial description of the relative motion between deputy satellites. To obtain the full description, we need to consider all $N-1$ linearly independent such assignments of coefficients, thus giving the total of $3(N-1)$ scalar equations. The simplest example of the set with suitable combinations of coefficients is
\[
 \{-1,1,0,\ldots,0\},\{0,-1,1,0,\ldots,0\},\ldots,\{0,\ldots,-1,1\}.
\]
Physically it means that we use relative displacements $\Delta {\mathbf r}_{i+1}- \Delta {\mathbf r}_i (i=1,N-1)$ of consecutive objects to describe the dynamics of the subsystem composed of the deputy satellites. One more opportunity is provided by the set
\[
\left\{1-\frac{1}{N},-\frac{1}{N},\ldots,-\frac{1}{N}\right\},\ldots,
\left\{-\frac{1}{N},\ldots,1-\frac{1}{N},-\frac{1}{N}\right\},
\]
characterizing relative displacements of the deputy satellites $D_1, \ldots,D_{(N-1)}$ with respect to the center of mass of all deputy satellites.
\end{enumerate}

Clearly, the collection of these three groups of equations is equivalent to the original system of $3(N+1)$ scalar equations for the main and deputy satellites. Let us deal with the three groups one by one.

\subsection{System CoM motion equations}

The motion of the system's center of mass is described by free HCW equations \cite{C1960}. It is well-known that these equations describe oscillations with frequency ${\omega }_0$ in the $y$ component, and oscillations with frequency ${\omega }_0$ combined with a linear drift in the orbital $xz$ plane.

\subsection{ Relative motion between the main satellite and the deputy satellites' center of mass }

Let $\Delta x,\Delta y,\Delta z$ denote the scalar components of $\frac{1}{N}\sum^N_{i=1}{\Delta {{\mathbf r}}_i}-{\Delta {\mathbf r}}_C$. Then using relation~\eqref{eq:len}

\begin{align}
\Delta \ddot{x} &= 2\omega_0 \Delta \dot{z} - \lambda_* \frac{k}{m_r}\Delta x=2\omega_0\Delta \dot{z}-3\omega_0^2\Delta x, \notag \\
\Delta \ddot{y} &=-\omega_0^2 \Delta y - \lambda_* \frac{k}{m_r} \Delta y = -4\omega_0^2\Delta y, \\
\Delta \ddot{z} &=-2\omega_0\Delta \dot{x} + 3\omega_0^2\Delta z - \frac{k }{m_r}\Delta z - \frac{b}{m_r}\Delta \dot{z}. \notag
\end{align}

The equation for $\Delta y$ is independent from the equations for $\Delta x,\Delta z$ and describes harmonic oscillations with the frequency $2{\omega }_0.$ To analyze the remaining equations for $\Delta x$, $\Delta z$ we consider the corresponding first-order system

\begin{equation} \label{eq:s1}
\frac{d}{dt}
\left(
\begin{array}{c}
\Delta x \\
\Delta \dot{x} \\
\Delta z \\
\Delta \dot{z}
\end{array}
\right)=\left(
\begin{array}{cccc}
0 & 1 & 0 & 0 \\
-3{\omega }^2_0 & 0 & 0 & 2{\omega }_0 \\
0 & 0 & 0 & 1 \\
0 & -2{\omega }_0 & 3{\omega }^2_0-\frac{k}{m_r} & -\frac{b}{m_r}
\end{array}
\right)
\left(
\begin{array}{c}
\Delta x \\
\Delta \dot{x} \\
\Delta z \\
\Delta \dot{z}
\end{array}
\right).
\end{equation}

Let us denote by $A$ the matrix on the right-hand side of the system (\ref{eq:s1}). To examine the
stability property of this system, we write down the characteristic equation
$$
\operatorname{det}(A-\rho I_4)=\rho^4 + \frac{b}{m_r}\rho^4 + \left(\frac{k}{m_r}+4\omega_0^2\right)\rho^2 +
\frac{3b\omega_0^2}{m_r}\rho + 3\omega_0^2\left(\frac{k}{m_r}-3\omega_0^2\right)=0.
$$

The symbol $I_4$ is used here for the identity matrix of the fourth order.
Applying the Routh-Hurwitz stability criterion \cite{G1959}, it is not difficult to establish
that all eigenvalues $\rho$ of $A$ belong to the left half-plane $Re \rho < 0$ iff the condition
\eqref{eq:stab1} is satisfied and $b>0$.

Summarizing, under the stability assumption \eqref{eq:stab1} the center of mass of deputy satellites oscillates with respect to the main satellite; in the linear approximation the energy dissipation in tethers does not affect the oscillations along
the $y$ axis.

\subsection{ Relative motion between deputy satellites }

Let $\Delta x,\Delta y,\Delta z$ denote the scalar components of $\sum^N_{i=1}{\sigma_i\Delta {{\mathbf r}}_i}$ with some assignment of coefficients $\sigma_i$ subject to $\sum^N_{i=1}{\sigma_i}=0$. Then

\begin{align}
\Delta \ddot{x}&=2\omega_0\Delta \dot{z} - \lambda_*\frac{k }{m_D}\Delta x, \notag \\
\Delta \ddot{y}&=-\omega_0^2\Delta y - \lambda_*\frac{k}{m_D}\Delta y, \\
\Delta \ddot{z}&=-2\omega_0\Delta \dot{x} + 3\omega_0^2\Delta z - \frac{k}{m_D}\Delta z - \frac{b}{m_D}\Delta \dot{z}.\notag
\end{align}

The equation for the $y$ component describes oscillations with frequency
\begin{equation}
\omega_y=\sqrt{\omega_0^2+\lambda_*\frac{k}{m_D}} = \omega_0\sqrt{\frac{4m_C+Nm_D}{m_C+Nm_D}}.
\end{equation}

The remaining equations for $\Delta x,\ \Delta z$ can be written in the first-order form as
\begin{equation} \label{eq:s2}
\frac{d}{dt}\left( \begin{array}{c}
\Delta x \\
\Delta \dot{x} \\
\Delta z \\
\Delta \dot{z} \end{array}
\right)=\left( \begin{array}{cccc}
0 & 1 & 0 & 0 \\
-\lambda_*\frac{k}{m_D} & 0 & 0 & 2{\omega }_0 \\
0 & 0 & 0 & 1 \\
0 & -2{\omega }_0 & 3{\omega }^2_0-\frac{k}{m_D} & -\frac{b}{m_D} \end{array}
\right)\left( \begin{array}{c}
\Delta x \\
\Delta \dot{x} \\
\Delta z \\
\Delta \dot{z} \end{array}
\right).
\end{equation}

The matrix $A_1$ on the right-hand side of \eqref{eq:s2} has the characteristic polynomial
\[
\operatorname{det}(A_1-\rho I_4) =\rho^4 + \frac{b}{m_D}\rho^3 + \left(\frac{(\lambda_*+1)k}{m_D} + \omega_0^2\right)\rho^2+\frac{\lambda_* k b}{m^2_D}\rho + \frac{\lambda_* k}{m_D}\left(\frac{k}{m_D}-3\omega_0^2\right).
\]

Applying again Routh-Hurwitz stability criterion, we establish that the system \eqref{eq:s2} is asymptotically
stable iff $b>0$ and
\begin{equation}\label{eq:stab2}
k\ge {3\omega_0^2m}_D.
\end{equation}

If the condition \eqref{eq:stab2} is fulfilled, but the dissipation is absent ($b=0$), then the spectrum of $A_1$ is purely
imaginary. If the condition is violated, matrix $A_1$ has a positive eigenvalue. Note that the stability condition \eqref{eq:stab2} is stronger than the earlier condition \eqref{eq:stab1}.

Assuming $b=0$, in the limit of large rigidity $k$ the linearized dynamics of $\Delta x$, $\Delta z$ approximately decouples into independent oscillations of $\Delta x$ and $\Delta z$ with frequencies
\begin{equation} \label{eq:omegaXZ}
\begin{aligned}
\omega_x &= \sqrt{\frac{\lambda_* k }{m_D}}\left(1-2\frac{m_D\omega_0^2}{k} + O\left(k^{-2}\right)\right) =\\
&=\sqrt{\frac{3m_C}{Nm_D+m_C}}\left(1-2\frac{m_D\omega_0^2}{k} + O\left(k^{-2}\right)\right)\omega_0,\\
\omega_z &= \sqrt{\frac{k}{m_D}}\left(1 + \frac{m_D\omega_0^2}{2k}  +O\left(k^{-2}\right)\right).
\end{aligned}
\end{equation}

\subsection{Lyapunov function}

Our solution of the linearized dynamics implies, in particular, that in the linear approximation of tether tension the vertical equilibrium is stable. This conclusion can also be shown in a stronger sense -- without linearization of tether tension -- by directly providing a Lyapunov function. Specifically, let $\mathbf{r}' = \mathbf{r} - \mathbf{r}_{CoM}$, $\dot{\mathbf{r}}'=\dot{\mathbf{r}} -\dot{\mathbf{r}}_{CoM}$ denote position and velocity relative to the system's center of mass. Consider the energy $E=\mathcal{T}+V$ of relative motion, where

\[\mathcal{T}=\frac{m_D}{2}\sum^N_{i=1}{\dot{{\mathbf r}'}_i}^2+\frac{m_C}{2}{{\dot{{\mathbf r}'}}_C}^2\]
and

\begin{align*}
V=&\frac{m_D{\omega }^2_0}{2}\sum^N_{i=1}{({{y'}^2_i-3z'}^2_i)}+\frac{m_C{\omega }^2_0}{2}\left({y'}^2_C-{3z'}^2_C\right)+\\
&+\frac{k}{2}\sum^N_{i=1}{\mathbf 1}_{(\left|{{\mathbf r}'}_C -{{\mathbf r}'}_{i}\right|>l_0)}
\cdot {\left(\left|{{\mathbf r}'}_{C}-{{\mathbf r}'}_i\right|-l_0\right)}^2.
\end{align*}

If $b=0,$ then $\dot{E}=0$, otherwise $\dot{E}\le 0$. Taylor expansion of $V$ near the equilibrium yields

\begin{align}
V&=V_{\rm equilibrium}+\frac{m_D{\omega_0}^2}{2}\sum^N_{i=1}{\left({\Delta y'_i}^2-3{\Delta z'_i}^2 \right)}+\frac{m_C{\omega_0}^2}{2}\left({\Delta y'_C}^2-{3\Delta z'_C}^2\right) + \notag \\
&+\frac{k}{2}\sum^N_{i=1}\left\{\left(\Delta {z'_i}-\Delta {z'_C}\right)^2+\lambda_*\left[{\left(\Delta {x'}_i-
\Delta {x'}_C\right)}^2+{\left(\Delta {x'}_i-\Delta {x'}_C\right)}^2\right]\right\} \\
&+ O({{|\Delta\mathbf{r}}'|}^3) \notag
\end{align}

Taking into account the identity
\[
\Delta{\mathbf{r}'}_C=-\frac{m_D}{m_C}\sum^N_{i=1}{\Delta{{\mathbf{r}}'}_i},
\]
one can check that this form is positive definite exactly if stability conditions \eqref{eq:stab1} and \eqref{eq:stab2} hold.

\section{Motion of deputy satellites along Lissajous curves}

Having described oscillations of the ``hub-and-spoke'' system near the equilibrium, we consider now the
possibility of the satellites moving in such a way that the system elements (satellites and tethers) never collide. Given our assumption of small deviations from the vertical equilibrium, we formulate this as the requirement that the projections $(x_i,y_i)$ of deputy satellites to the $xy$ plane never come close to each other.

Results of Section 3.3 imply that the position of a deputy satellite relative to the center of mass of all deputy satellites,
\[{{\mathbf r}'}_i={\mathbf r}_i-\frac{1}{N}\sum^N_{k=1}{{\mathbf r}_k},\]
 oscillates with frequency
\begin{equation}\label{eq:omegaY}
\omega_y=\omega_0\sqrt{\frac{{4m}_C+Nm_D}{m_C+Nm_D}}
\end{equation}
in the $y$ direction and, for sufficiently rigid tethers, with frequency
\begin{equation}\label{eq:omegaX}
\omega_x \approx \omega_0\sqrt{\frac{3m_C}{Nm_D+m_C}}
\end{equation}
in the $x$ direction (with a more accurate value given by \eqref{eq:omegaXZ}). Since, as shown in Section 3.2, the center of mass of deputy satellites $\frac{1}{N}\sum^N_{k=1}{{{\mathbf r}}_k}$
performs independent oscillations, we may ignore these latter and assume without loss of generality that
\begin{equation}\label{eq:eq4}
\frac{1}{N}\sum^N_{k=1}{x_k}=\frac{1}{N}\sum^N_{k=1}{y_k}\equiv 0
\end{equation}
at all times.

By linearity, for each pair $i,j$ of deputy satellites their relative position ${\mathbf r}_i-{\mathbf r}_j$ also oscillates with the same frequencies $\omega_x$, $\omega_y$. If the frequencies are incommensurate, i.e.,  $\omega_x/\omega_y$ is irrational, then the trajectory of the oscillation is aperiodic and it comes arbitrarily close to the origin, i.e. the two satellites come arbitrarily close to each other. We are thus naturally led to consider commensurate oscillations:
\begin{equation}\label{eq:eq5}
\frac{\omega_x}{\omega_y}=\frac{p}{q},
\end{equation}
where $p, q$ are co-prime natural numbers. In this case the two oscillations have a common period of
\[T_L=\frac{2\pi p}{\omega_x}=\frac{2\pi q}{\omega_y}.\]

It is convenient to introduce the non-dimensional time
\[\tau =\frac{t}{T_L}.\]

The $xy$-motion of a single deputy satellite can then be written as
\begin{equation}\label{eq:eq_single_ds}
x=x_0 \sin  \left(2\pi p\tau + \varphi_x\right),\ \ y=y_0\sin  (2\pi q\tau + \varphi_y),
\end{equation}
with some initial phases $\varphi_x,\ \varphi_y$, and has a period 1. The trajectory of this motion is known as a Lissajous curve \cite[Sec. 25]{A1978}

Substituting expressions \eqref{eq:omegaY}, \eqref{eq:omegaX} for $\omega_x$, $\omega_y$ into formula \eqref{eq:eq5}, we obtain the following relation between the frequency ratio and the satellite mass ratio:
\[
\frac{Nm_D}{m_C}=\frac{3q^2}{p^2}-4.
\]

In particular, positivity of the left-hand side entails
\begin{equation}\label{eq:eq7}
\frac{p}{q} < \frac{\sqrt{3}}{2}.
\end{equation}

We remark in passing that this condition excludes the usual elliptic (or circular) oscillations corresponding to $p=q=1$.

The frequencies $\omega_x,\ \omega_y$ can then be expressed in terms of $p$ and $q$:
\[\omega_x=\frac{p\omega_0}{\sqrt{q^2-p^2}}, \qquad \omega_y=\frac{q\omega_0}{\sqrt{q^2-p^2}}.\]

\subsection{Balanced formations avoiding collisions}

We seek now formations of deputy satellites moving according to \eqref{eq:eq_single_ds} subject to the following conditions:

\begin{itemize}
\item[{\bf A.}]  The arrangement of deputy satellites must satisfy the balance condition (\ref{eq:eq4});

\item[{\bf B.}]  The satellites must never collide, i.e.
\[
(x_i(t),y_i(t))\ne (x_j(t),y_j(t))
\]
  for all $t$ and $i\ne j$;

\item[{\bf C.}]  Optionally, we may wish to ensure that
\[(x_i(t),y_i(t))\ne (0,0)\]
for all $i$ and $t$ -- in this case an additional satellite can be added at the center of the ``hub-and-spoke'' system without collisions with this system.
\end{itemize}

We consider two types of uniform arrangement of deputy satellites. Type I is a uniform arrangement of  $N$ satellites along a single Lissajous curve: the position of the $i$'th deputy satellite is given by
\[x_i(\tau )=x_0 \sin  \left[2\pi p\left(\tau +\frac{i}{N}\right)+{\varphi }_x\right],
\qquad y_i(\tau )=y_0 \sin  \left[2\pi q\left(\tau +\frac{i}{N}\right)+{\varphi }_y\right],\ \]
where the phases $\varphi_x,\ \varphi_y$ are the same for all satellites.

Type II is a uniform arrangement of $N$ satellites along several Lissajous curves: the position of the $i$'th deputy satellite is given by
\[x_i(\tau )=x_0\sin  \left[2\pi \left( p\tau +\frac{i}{N}\right)+{\varphi }_x\right], \qquad y_i(\tau )=y_0 \sin  \left[2\pi \left( q\tau +\frac{i}{N}\right)+{\varphi }_y\right],\]
where the phases $\varphi_x,\ \varphi_y$ are the same for all satellites.

The following proposition summarizes properties of such formations with respect to the above three conditions.

\begin{prop} Let $N=2,3,\dots $ Denote $\varphi_0=\frac{q\varphi_x-p\varphi_y}\pi$.
\begin{itemize}
\item[a)] For a Type I formation, the balance condition A is fulfilled if and only if neither $p$ nor $q$ is divisible by $N$. For a Type II formation, the balance condition A is fulfilled for all $N$.

\item[b)]  For a Type I formation, the no-collision condition B is fulfilled iff the number $\varphi_0+(p-q)/2$ is not an integer and $N$ is co-prime with $p$ and $q$.

\item[c)]  For a Type II formation, in case $N\ge 3$ the no-collision condition B is fulfilled if and only if $\left(\varphi_0+(p-q)/2\right)N$ is not divisible by the greatest common divisor of $N$ and $q-p$. In case $N=2$ condition B is fulfilled iff $\varphi_0$ is not an integer.

\item[d)]  Lissajous curve \eqref{eq:eq_single_ds} goes through the origin $(0,0)$ iff $\varphi_0$ is an integer. It follows that a Type I formation satisfies condition C iff $\varphi_0$ is not an integer. A Type II formation satisfies condition C iff $\varphi_0$ is not of the form $a+2b(q-p)/N$ with integer $a,b$.
\end{itemize}
\end{prop}

The proof of Proposition 1 is given in Appendix \ref{sec:A}. Note that, given $p,q$, admissible formations (satisfying all three conditions A-C) exist for all  $N=2,3,\dots $ in case of Type II, but not for all $N$ in case of Type I.

In Table 1 we list parameters of admissible formations for all $p,q\le 4$ subject to condition (\ref{eq:eq7}).

\begin{table}[hb]
\begin{center}
\begin{tabular}{cll}
\toprule
$p/q(=\omega_x/\omega_y)$ & $Nm_D/m_C$ & Admissible $N$ for Type I \\
\midrule
$1/2$ & \hspace{0.5cm}8    & \hspace{1cm}3, 5, 7, \dots  \\
$1/3$ & \hspace{0.5cm}23   & \hspace{1cm}2, 4, 5, \dots   \\
$2/3$ & \hspace{0.5cm}11/4 & \hspace{1cm}5, 7, 11, \dots  \\
$1/4$ & \hspace{0.5cm}44   & \hspace{1cm}3, 5, 7, \dots  \\
$3/4$ & \hspace{0.5cm}4/3  & \hspace{1cm}5, 7, 11, \dots \\
\bottomrule
\end{tabular}
\caption{Admissible formations for small values of $p,q$.}
\end{center}
\end{table}

\begin{center}
\begin{figure}[htb]
		\includegraphics[width=0.65\textwidth, keepaspectratio, angle=-90, clip, trim=0mm 10mm 10mm 5mm]{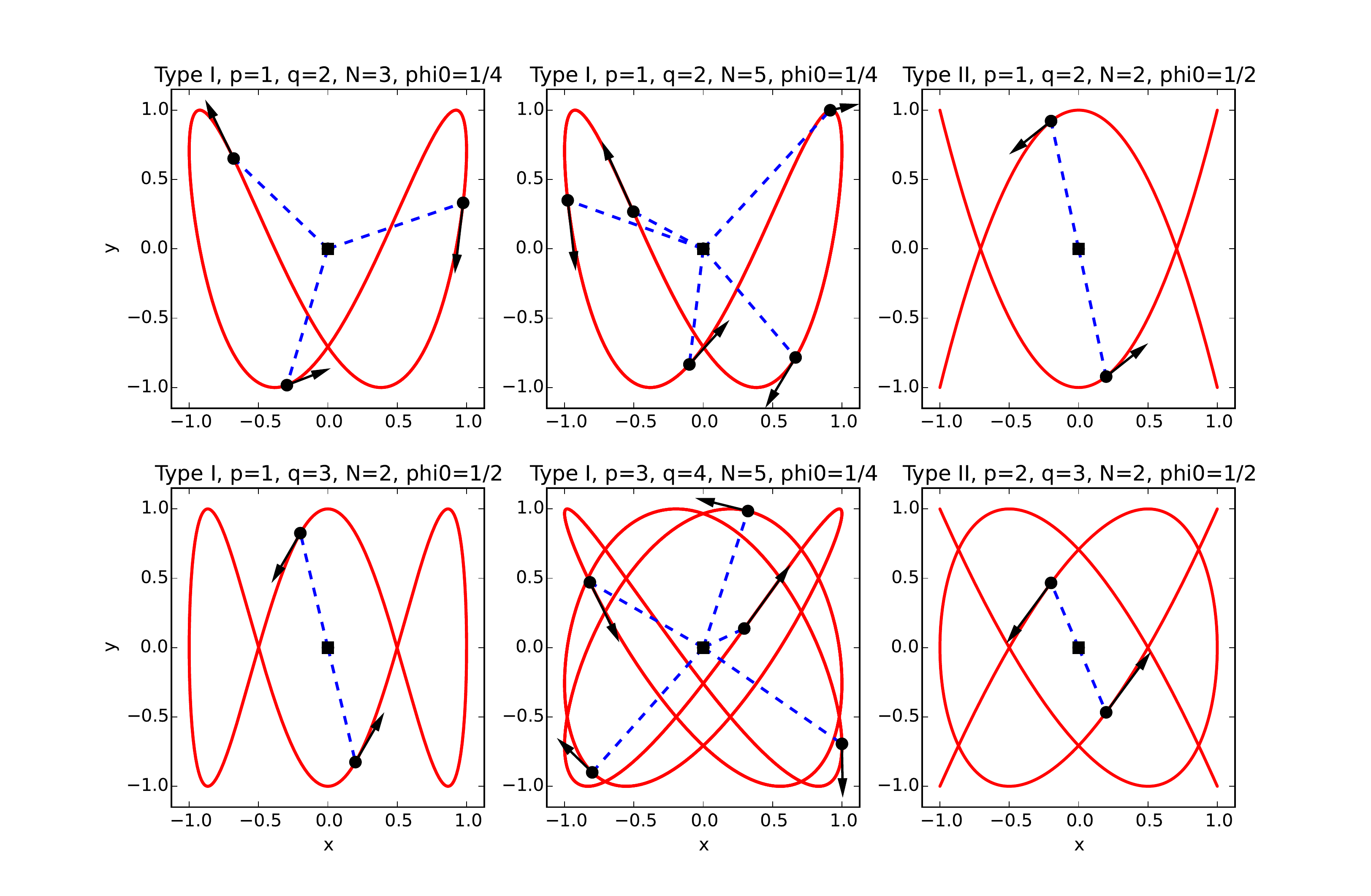}		
    \caption{Examples of formations of Types I and II. }
    \label{fig:fig2}
\end{figure}
\end{center}

In Fig. \ref{fig:fig2} we show several examples of formations of Types I and II satisfying conditions A,B,C with $p$, $q$ subject to condition \eqref{eq:eq7}.

\subsection{Entanglement}

Practical implementation of the introduced formations requires to resolve another issue. The tethers' ends are not attached to the main satellite at exactly the same point. Even if the tethers are connected very close to each other, they still have nonzero thickness. Relative motion of satellites in formations of Type I and Type II may cause not only contacts between the tethers, but also tethers entanglement.

\begin{figure}[htb]
\centering
		\includegraphics[width=0.95\textwidth, keepaspectratio]{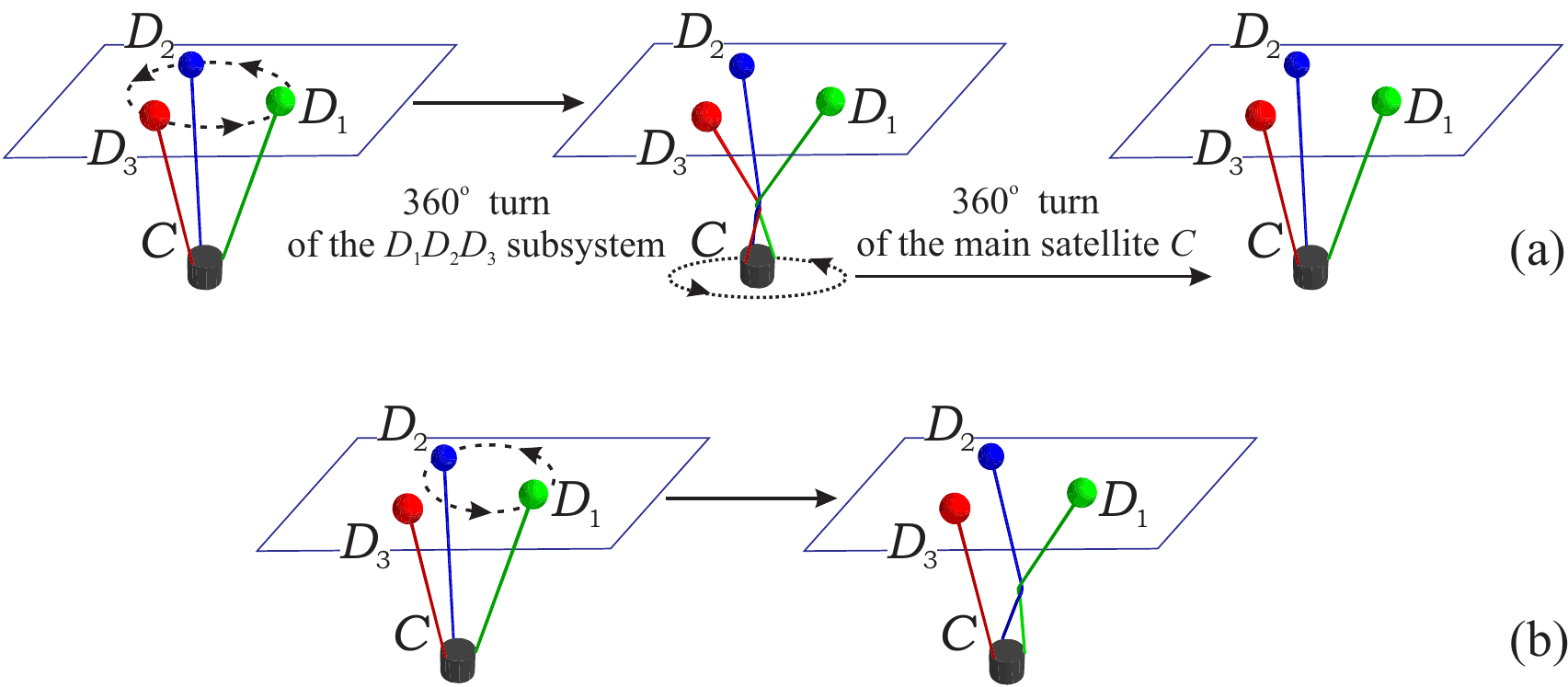}		
    \caption{Examples of weak and strong entanglement. In the case of weak entanglement (a) tethers can be disentangled by rotation of the main satellite. Tethers in strong entanglement (b) can not be disentangled by rotation of the main satellite.}
    \label{fig:EntAB}	
\end{figure}

We will distinguish two kinds of entanglement: one that can be canceled out by rotating the main satellite about the $z$ axis as shown in Fig.~\ref{fig:EntAB}a, and one that can not be eliminated by such rotations (Fig.~\ref{fig:EntAB}b). These two kinds will be referred to as \textit{weak} and \textit{strong} entanglement, respectively.

We can define entanglement rigorously by applying simple topological concepts \cite{A1983} to the dynamics of the deputy satellites. Consider the motion of a deputy satellites formation as a continuous map 
\[
{{\mathcal M}}_0\left(\tau\right) = \left(\left( \begin{array}{c}
 x_1(\tau) \\
 y_1(\tau) \end{array}
\right),\left( \begin{array}{c}
 x_2(\tau) \\
 y_2(\tau) \end{array}
\right) ,\dots \left(\begin{array}{c}
 x_N(\tau) \\
 y_N(\tau) \end{array}
\right)\right),\ \ \tau \ \in [0,1].
\]
By periodicity, $\mathcal{M}_0\left(0\right)=\mathcal{M}_0\left({\rm 1}\right)$. Entanglement means the impossibility to ``straighten'' this map subject to topological constraints. Precisely, let us say that ${{\mathcal M}}_0$ is \textit{strongly homotopic} to the constant map
\[\mathcal{M}_1\left(\tau\right) =\left(\left( \begin{array}{c}
x_1(0) \\
y_1(0) \end{array}
\right),\left( \begin{array}{c}
x_2(0) \\
y_2(0) \end{array}
\right),\dots \left( \begin{array}{c}
x_N(0) \\
y_N(0) \end{array}
\right)\right),\ \ \tau \ \in \left[0,1\right],\]
if there is a continuous family of maps $\mathcal{M}_{\alpha }\left(\tau \right),\tau {\rm \ }\in \left[{\rm 0,1}\right],\alpha {\rm \ }\in \left[{\rm 0,1}\right],\ $that deforms $\mathcal{M}_0$ into $\mathcal{M}_1$ without collisions and subject to the periodicity constraint
\begin{equation}\label{eq:eq8}
\mathcal{M}_{\alpha }\left(0\right)=\mathcal{M}_{\alpha }\left({\rm 1}\right),\ \ \alpha \ \in [0,1].
\end{equation}
If $\mathcal{M}_0$ is not strongly homotopic to $\mathcal{M}_{\rm 1}$, then we call $\mathcal{M}_0$ \textit{weakly entangled}.

In order to define strong entanglement, we say that $\mathcal {M}_0$ is \textit{weakly homotopic} to the constant map $\mathcal{M}_{\rm 1}$ if the constraint \eqref{eq:eq8} is replaced by the relaxed constraint
\[\mathcal{M}_{\alpha }\left(0\right)=R_{\varphi(\alpha)}\mathcal{M}_{\alpha }\left({\rm 1}\right),\ \ \alpha \ \in \left[0,1\right],\ \]
where $R_{\varphi (\alpha )}$ is the rotation of the satellite positions in the $xy$ plane by continuously varying  angles $\varphi (\alpha {\rm )}$ common to all deputy satellites, for some choice of the continuous function $\varphi (\alpha {\rm )}$. We then say that $\mathcal{M}_0$ is \textit{strongly entangled} if $\mathcal{M}_0$ is not weakly homotopic to $\mathcal{M}_{\rm 1}$.

The above definitions do not involve the main satellite. We assume that the initial state is non-entangled, and the geometry of the tether attachment points to the main satellite is reflected in the initial positions of the deputy satellites as shown in Fig.~\ref{fig:EntAB}.

An obvious example of a weak, but not strong entanglement results from a circular motion as shown in Fig.~\ref{fig:EntAB}a. Also note that a motion of $N=2$ deputy satellites is never strongly entangled but may be weakly entangled.

By considering a few simple examples, it is easy to see that entanglement does take place in some formations of Type I and Type II and does not in others. For instance, there is a strong braid-like entanglement for the Type I formation with $p=1$, $q=2$, and  $N=3$ (top left plot in Fig.~\ref{fig:fig2}), weak entanglement for the Type I formation with $p = 1$, $q = 3$, and $N = 2$ (bottom left plot in Fig.~\ref{fig:fig2}), and no entanglement for the Type II formations with $N = 2$ and $p=1$, $q=2$ or $p=2$, $q=3$, and  $N=2$ (top and bottom right plots in Fig.~\ref{fig:fig2}).

Complete analysis of entanglement in our formations appears to be relatively complex mathematically, hence we will not attempt it here. We will, however, state a simple proposition involving pairwise relations between satellites.  Namely, for any pair $(i,j)$ of deputy satellites consider the \textit{winding number} $w_{i,j}\ $defined as the number of turns that the $j$'th satellite makes about the $i$'th satellite in the $xy$-plane over the period:

\[w_{i,j}=\frac{1}{2\pi }{\left.{\arg  \left( \begin{array}{c}
x_j(\tau )-x_i(\tau ) \\
y_j(\tau )-y_i(\tau ) \end{array}
\right)\ }\right|}^1_{\tau =0}\]
Clearly, $w_{i,j}=w_{j,i}$.The winding numbers are obviously invariant under a strong homotopy. Under a weak homotopy, they are incremented by an amount common to all pairs and equal to the number of full turns  $\frac{1}{2\pi }{\left.\varphi (\alpha {\rm )}\right|}^1_{\alpha =0}$. Since winding is absent in the straightened map $\mathcal{M}_{\rm 1}$, the numbers $w_{i,j}$ can serve to establish sufficient (but not necessary) conditions of entanglement.

\begin{prop}\mbox{}\hfill
\begin{itemize}
\item[a)] Consider a formation of Type I or Type II moving without collisions. If at least one of $p$, $q$ is even, then $w_{i,j}=0$ for all pairs. Otherwise, $w_{i,j}=\pm 1$ for all pairs.Consequently, if both $p$, $q$ are odd, then there is at least a weak entanglement.
\item[b)]  Moreover, consider a formation of Type I and suppose that both $p$, $q$ are odd so that $w_{i,j}=\pm 1$ by a). Then both values  $+1$ and $-1$ are encountered among the winding numbers $w_{i,j}$ iff neither of the numbers $q-p,\ q+p$ is divisible by $2N$. As a consequence, if neither of $q-p,\ q+p$ is divisible by $2N$, then the entanglement is not only weak, but also strong.
\end{itemize}
\end{prop}

Proposition 2 is proved in Appendix \ref{sec:B}. It leads to an expected result if applied to a circular or elliptic motion ($p=q=1$): part a) confirms a weak entanglement, and part b) ensures that all winding numbers are equal since $p-q=0$, so the sufficient condition for a strong entanglement is not fulfilled.

Note also that for odd $p,\ q$ one of the numbers $q-p,\ q+p$ is always divisible by 4, so statement b) agrees with our earlier remark that there can be no strong entanglement for $N=2$.

\subsection{Second order perturbation theory}

In contrast to the general stability of the ``hub-and-spoke'' system established in Section 3.4 in the sense of smallness of deviation from the vertical equilibrium, we do not expect the motion of deputy satellites along Lissajous curves to be stable. This motion is a subtle phenomenon which, in particular, is strongly affected by nonlinearities and can be maintained without additional control action only at relatively small oscillation amplitudes.

In this section we examine the nonlinearity effects by deriving second order corrections to the evolution equations in case
$k/3\omega_0^2 m_D\gg 1$ when one can neglect the tethers' extensibility and put $l_0=l_*$.
If we assume that the center of mass of the system rests at the origin, then the configuration can be parametrized by the $2N$ coordinates $x_i, y_i,\ i=1,..,N,$ of the auxiliary satellites. Specifically, the remaining coordinates are given by
\begin{align*}
x_C &= -\frac{m_D}{m_C}\sum^N_{i=1}{x_i},\ y_C = -\frac{m_D}{m_C}\sum^N_{i=1}{y_i}, \\
z_C &=\frac{m_r}{m_C}\sum^N_{i=1}{\sqrt{l^2-{\left(x_i-x_C\right)}^2-{\left(y_i-y_C\right)}^2}}
\approx \\
&\approx \frac{m_r}{m_C}\left\{Nl_0-\frac{1}{2l_0}\sum^N_{i=1}{\left[{\left(x_i-x_C\right)}^2+{\left(y_i-y_C\right)}^2\ \right]}\right\},\\
z_i&=z_C-\sqrt{l^2-{\left(x_i-x_C\right)}^2-{\left(y_i-y_C\right)}^2} \approx -\frac{m_r}{m_D}l_0 + \\
&+\frac{1}{2l_0}\left\{\left[{\left(x_i-x_C\right)}^2+{\left(y_i-y_C\right)}^2\ \right] - \frac{m_r}{m_C}\sum^N_{k=1}{\left[{\left(x_k-x_C\right)}^2+{\left(y_k-y_C\right)}^2\ \right]}\right\}.
\end{align*}
Here and below we consistently keep terms only up to second order in $x_i,\ y_i, {\dot{x}}_i,\ {\dot{y}}_i,\ {\ddot{x}}_i,\ {\ddot{y}}_i$.

The $x,y$ components $T_{i,x},\ T_{i,y}$ of the tension acting on the $i$'th deputy satellite can be expressed through the $z$ component by
\[T_{i,x}=\frac{x_i-x_C}{z_i-z_C}T_{i,z}\approx -\frac{x_i-x_C}{l_0}T_{i,z},\ \ T_{i,y}=\frac{y_i-y_C}{z_i-z_C}T_{i,z}\approx -\frac{y_i-y_C}{l_0}T_{i,z}.\]

The $z$ component $T_{i,z}$, in turn, is found from the corresponding HCW equation:
\[\frac{T_{i,z}}{m}\approx 2{\omega }_0{\dot{x}}_i-3{\omega }^2_0z_i\approx 2{\omega }_0{\dot{x}}_i+3{\omega }^2_0\frac{m_r}{m_D}l_0.\]

The remaining two HCW equations then yield
\[{\ddot{x}}_i\approx 2{\omega }_0{\dot{z}}_i-\frac{x_i-x_C}{l_0}\left(2{\omega }_0{\dot{x}}_i+3{\omega }^2_0\frac{m_r}{m_D}l_0\right),\ \]
\[{\ddot{y}}_i\approx -{\omega }^2_0y_i-\frac{y_i-y_C}{l_0}\left(2{\omega }_0{\dot{x}}_i+3{\omega }^2_0\frac{m_r}{m_D}l_0\right).\]

Differentiating $z_i$ and retaining our notation for $\omega_x,\ \omega_y$ from (\ref{eq:omegaY})-(\ref{eq:omegaX}), we finally obtain
\begin{align*}
\ddot{x}_i+\omega_x^2\left(x_i-x_C\right) \approx & \frac{2\omega_0}{l_0}\bigg(
-{\dot{x}}_C\left(x_i-x_C\right)+\left({\dot{y}}_i-{\dot{y}}_C\right)\left(y_i-y_C\right)-\\
&-\frac{m_r}{m_C}\sum^N_{j=1}\big[\left(\dot{x}_j-\dot{x}_C\right)\left(x_j-x_C\right)
+\left(\dot{y}_j-\dot{y}_C\right)\left(y_j-y_C\right)\big]
\bigg),\\
\ddot{y}_i+\omega_y^2 y_i-3\omega_0^2\frac{m_r}{m_D}y_C\approx & -\frac{2\omega_0}{l_0}\dot{x}_i\left(y_i-y_C\right),
\end{align*}
where we have placed first order terms on the left and second order terms on the right.

The obtained equations can be used to find anharmonic corrections to a particular small harmonic oscillation. The usual procedure is to substitute the harmonic oscillation in the right-hand side and find the correction, to leading order, by solving the resulting non-homogeneous linear equation \cite{N1973}. If the right-hand side contains secular terms, i.e., those whose frequencies match some of the eigenfrequencies of the linear equation, then, additionally, the solution's frequencies need to be adjusted to prevent its non-physical growth.

In the case at hand we take the motion of a Type I or II formation with a small amplitude as a base harmonic oscillation that we denote $x_{i,0}\left(t\right),\ y_{i,0}\left(t\right)$. From the balance condition (\ref{eq:eq4}) we have $x_C\left(t\right)=y_C(t)\equiv 0$, so that the evolution equations simplify to
\begin{align*}
{\ddot{x}}_i+{\omega }^2_x\left(x_i-x_C\right) &\approx \frac{2{\omega }_0}{l_0}\left[{\dot{y}}_{i,0}y_{i,0}-\frac{m_r}{m_C}\sum^N_{j=1}{\left({\dot{x}}_{j,0}x_{j,0}+{\dot{y}}_{j,0}y_{j,0}\ \right)}\right],\\
\ddot{y}_i+\omega_y^2 y_i-3\omega_0^2\frac{m_r}{m_D}y_C &\approx -\frac{2{\omega }_0}{l_0}{\dot{x}}_{i,0}y_{i,0}.
\end{align*}

The linear terms describe oscillations with frequencies
\[{\omega }_x,\ {\omega }_y,\ {\omega }_{Cx}=\sqrt{3}\omega_0,{\omega }_{Cy}=2\omega_0,\]
where the latter two correspond to the motion of the center of mass of auxiliary satellites or, equivalently, to the motion of the main satellite. The terms ${\dot{x}}_{k,0}x_{k,0},\ \ {\dot{y}}_{k,0}y_{k,0}$, and ${\dot{x}}_{i,0}y_{i,0}$ on the right-hand side result in oscillations with frequencies $2{\omega }_x,2{\omega }_y$, and ${\omega }_y\pm {\omega }_x$, respectively.

Though the second order correction of motion is obviously present for each deputy satellite for any choice of system parameters, it is possible to choose parameters so as to make the second order correction completely vanish for the main satellite. The evolution equations for the main satellite are derived by adding up the equations for the deputy satellites:
\begin{align*}
-\frac{m_C}{m_D}\left({\ddot{x}}_C+{\omega }^2_{Cx}x_C\right)&\approx \frac{2{\omega }_0}{l_0}\left(\frac{m_C}{m_C + N m_D}\sum^N_{j=1}{{\dot{y}}_{j,0}y_{j,0}}-\frac{N m_D}{m_C+N m_D}\sum^N_{j=1}{{\dot{x}}_{j,0}x_{j,0}}\right),\\
-\frac{m_C}{m_D}\left({\ddot{y}}_C+{\omega }^2_{Cy}y_C\right)&\approx -\frac{2{\omega }_0}{l_0}\sum^N_{j=1}{{\dot{x}}_{j,0}y_{j,0}}.
\end{align*}

Note that the $Nm_D$ tends to be much larger than $m_C$ for a system satisfying our assumptions (see Table 1), so the motion of the main satellite tends to be generally affected by second order corrections much stronger than the motion of the deputy satellites. However, contributions from different $j$'s here may cancel out.

\begin{prop}  Consider a Type I formation such that neither of $2p,\ 2q,\ q\pm p$ is divisible by $N$. Then
\[
\sum^N_{k=1}{{\dot{x}}_{k,0}(t)x_{k,0}(t)}=\sum^N_{k=1}{{\dot{y}}_{k,0}(t)y_{k,0}(t)}=\sum^N_{k=1}{{\dot{x}}_{k,0}(t)y_{k,0}(t)}\equiv 0
\]
for all $t$, so that the right-hand sides of the above equations vanish identically.
\end{prop}

The proof is elementary, and we omit it. Examples of parameter sets fulfilling the proposition's hypothesis are $N=5$, $p=1$, $q=2$ and $N=5$, $p=3$, $q=4$. Our numerical experiments below confirm that in these cases the main satellite is indeed stable in contrast to the generic settings at the same oscillation amplitude of deputy satellites.

\section{Examples and numerical simulations}
\subsection{Setup and general observations}
In our numerical simulations we consider oscillations of the system with rigid tether at small angles and with equal amplitudes in $x$ and $y$:
\[
x_0=y_0=a=\varkappa_{\rm rad} l_0
\]

Here $a$ is the linear amplitude and $\varkappa_{\rm rad}$ is the corresponding angular amplitude expressed in radians; the same angular amplitude expressed in degrees is denoted $\varkappa_{\rm deg}$. In all our experiments
the dimensionless coefficient  $k/3\omega_0^2 m_D$  characterizing tether rigidity falls in the range $[3 \cdot 10^2,10^3]$
and $\varkappa_{\rm deg} \le 6^\circ$ (to justify the above linear relation between $a$ and $\varkappa_{\rm rad}$).
In order to numerically examine the stability of the Lissajous motion we introduce quantities characterizing relative deviations of the main and deputy satellites from their theoretical positions obtained in the linear approximation. Specifically, we consider the relative deviation of the main satellite's numerically computed position from the $z$ axis
\[
\delta_C(t)=\frac{1}{a}\sqrt{x^2_{C,{\rm num}}(t)+y^2_{C,{\rm num}}(t)},
\]
and the mean relative deviation of the deputy satellites' numerically found trajectories from the theoretical Lissajous curves:
\[
\delta_D(t)=\frac{1}{Na}\sum_{i=1}^{N}\sqrt{\left(x_{i,{\rm num}}(t)-x_{i,{\rm Liss}}(t)\right)^2+\left(y_{i,{\rm num}}(t)-y_{i,{\rm Liss}}(t)\right)^2}
\]
These quantities can be compared with the minimum distance between different deputy satellites on their theoretical trajectories:
\[
\delta_{\rm min}=\frac{1}{a}\min_{t,i \ne j}\sqrt{\left(x_{i,{\rm Liss}}(t)-x_{j,{\rm Liss}}(t)\right)^2+\left(y_{i,{\rm Liss}}(t)-y_{j,{\rm Liss}}(t)\right)^2}
\]
Here we neglect the difference in $z$ coordinates of the deputy satellites since for small amplitude oscillations this difference has higher order of smallness.

If each satellite remains at all times within a distance of $\delta_{\rm min}/2$ from its theoretical position on the Lissajous curve, then all satellites are guaranteed to avoide collisions. We will consider a slightly relaxed condition of stability

\begin{equation}\label{eq:collision}
\delta_D(t)<\frac{\delta_{\rm min}}{2}
\end{equation}
that constrains only the mean deviation $\delta_D(t)$ of the deputy satellites. The time interval during which this condition holds can be roughly considered  as a ``system stability interval''.

We restrict ourselves to the simplest ratio $\omega_x:\omega_y =1:2$ achieved at $p=1$, $q=2$. In this case, the theoretical mass ratio is
\[
\frac{Nm_D}{m_C}=8,
\]
and the full period $T_L$ of system oscillations is related to the orbital period $T_0=1/(2\pi\omega_0)$ by
\[
T_L=\sqrt{q^2-p^2}T_0={\sqrt{3}}T_0,
\]
so  that the non-dimensional time
\[
\tau = \frac{t}{\sqrt{3}T_0}.
\]

We consider the three different configurations of the system shown in the first row of Fig.~\ref{fig:fig2}. For each configuration, we perform simulations for $\varkappa_{\rm deg}=1^\circ$ and $\varkappa_{\rm deg}=3^\circ$ for 10 orbital periods. In our simulations, the central satellite moves along a geostationary orbit, and tether length $l_0=10000$ m.

\begin{figure}[htb]
\begin{center}
    \includegraphics[width=0.8\textwidth, keepaspectratio]{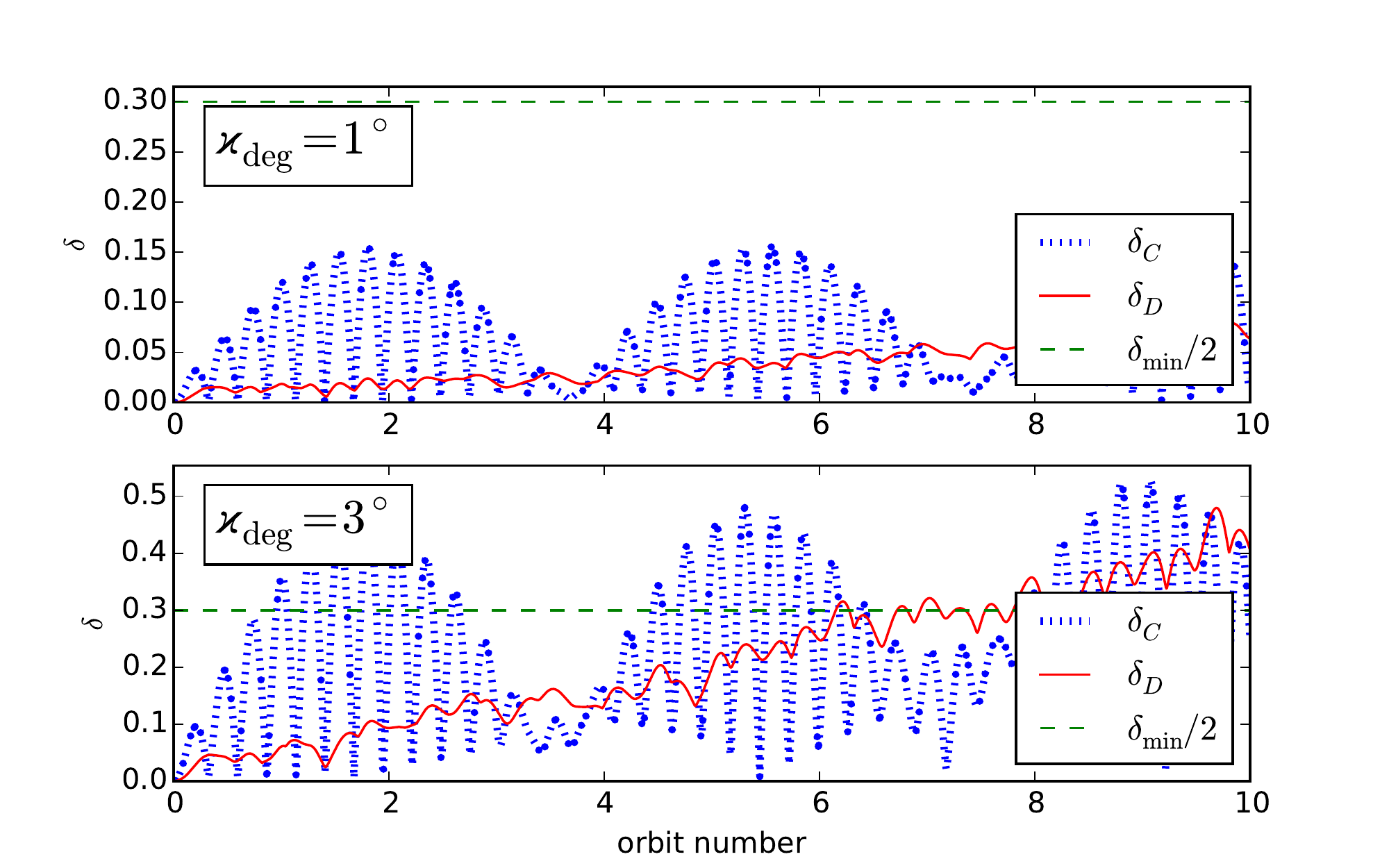}
\end{center}
    \caption{Numerically observed relative deviations of the satellites from the theoretical trajectories for the Type I formation with $N=3$ deputy satellites. In this case the  minimum relative distance between deputy satellites on the theoretical trajectories is $\delta_{\rm min}=0.60$. At $\varkappa_{\rm deg}=1^\circ$ stability condition \eqref{eq:collision} holds with a large margin for the whole simulation interval of 10 orbital periods, while at $\varkappa_{\rm deg}=3^\circ$ it breaks down after six orbital periods. }
    \label{fig:devTypeI_N3}
\end{figure}

The results for Type I formation at $N=3$ (Fig.~\ref{fig:devTypeI_N3})
show that the deviations of the main satellite are initially much larger than those of the deputy satellites. This is not surprising, since the main satellite is eight times lighter. However, deviations of the deputy satellites approximately linearly accumulate with time, and eventually catch up with those of the main satellite.
In case $\varkappa_{\rm deg}=3^\circ$ relative deviations are much larger than in case $\varkappa_{\rm deg}=1^\circ$: approximately three times larger for $\delta_C$ and six times larger for $\delta_D$. The second order perturbation theory in Section 4.3 suggests a linear dependence of relative deviations on $\varkappa_{\rm deg}$, but for $\varkappa_{\rm deg}=3^\circ$ deviations of the main satellite from the equilibrium position are already comparable to the oscillation amplitude, so this perturbation theory is not truly applicable here.
In case $\varkappa_{\rm deg}=1^\circ$ the stability condition \eqref{eq:collision} holds with a large margin for the whole simulation interval, while at $\varkappa_{\rm deg}=3^\circ$ it breaks down after six orbital periods.

\begin{figure}[htb]
\begin{center}
    \includegraphics[width=0.8\textwidth, keepaspectratio]{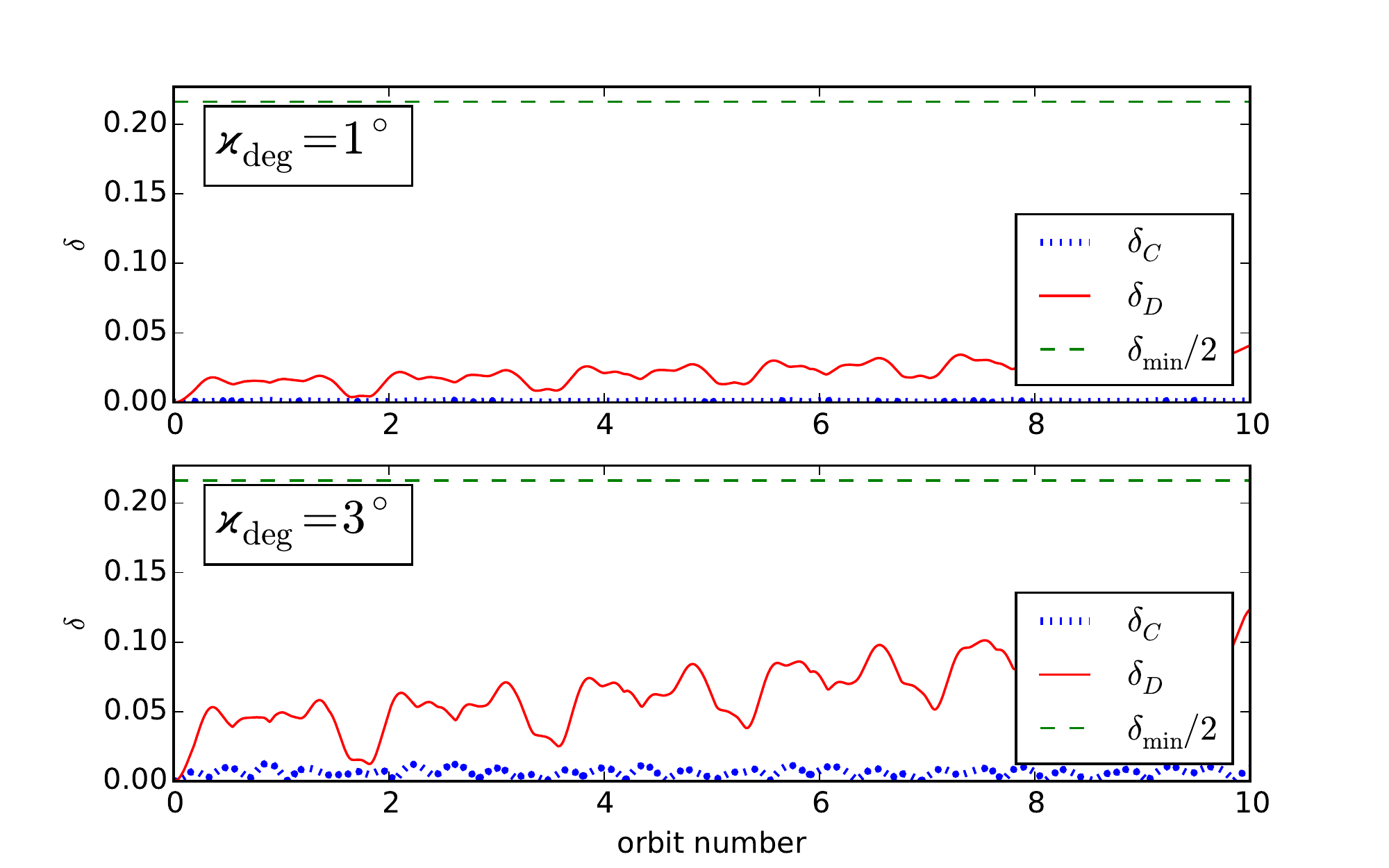}
\end{center}
    \caption{Deviations for Type I and $N=5$. In this case $\delta_{\rm min}=0.43$. Stability condition \eqref{eq:collision} holds with a large margin for the whole simulation interval of 10 orbital periods for both $\varkappa_{\rm deg}=1^\circ$ and $\varkappa_{\rm deg}=3^\circ$.}
    \label{fig:devTypeI_N5}
\end{figure}

The results for Type I formation at $N=5$ (Fig.~\ref{fig:devTypeI_N5}) are drastically different due to the cancellation of second order corrections of the main satellite's motion pointed out in Proposition 3. Not only is the main satellite almost immobile, but also the deviations of the deputy satellites are from two to four times smaller than in the previous case. In particular, the no-collision condition \eqref{eq:collision} holds at $\varkappa_{\rm deg}=3^\circ$ throughout the whole simulation interval.
Note also that in the $\varkappa_{\rm deg}=3^\circ$ case the deviations of deputy satellites are approximately three times as large as in the $\varkappa_{\rm deg}=1^\circ$ case, in good agreement with the second order perturbation theory which is now applicable since the deviations are small.

\begin{figure}[htb]
\begin{center}
		\includegraphics[width=0.8\textwidth, keepaspectratio]{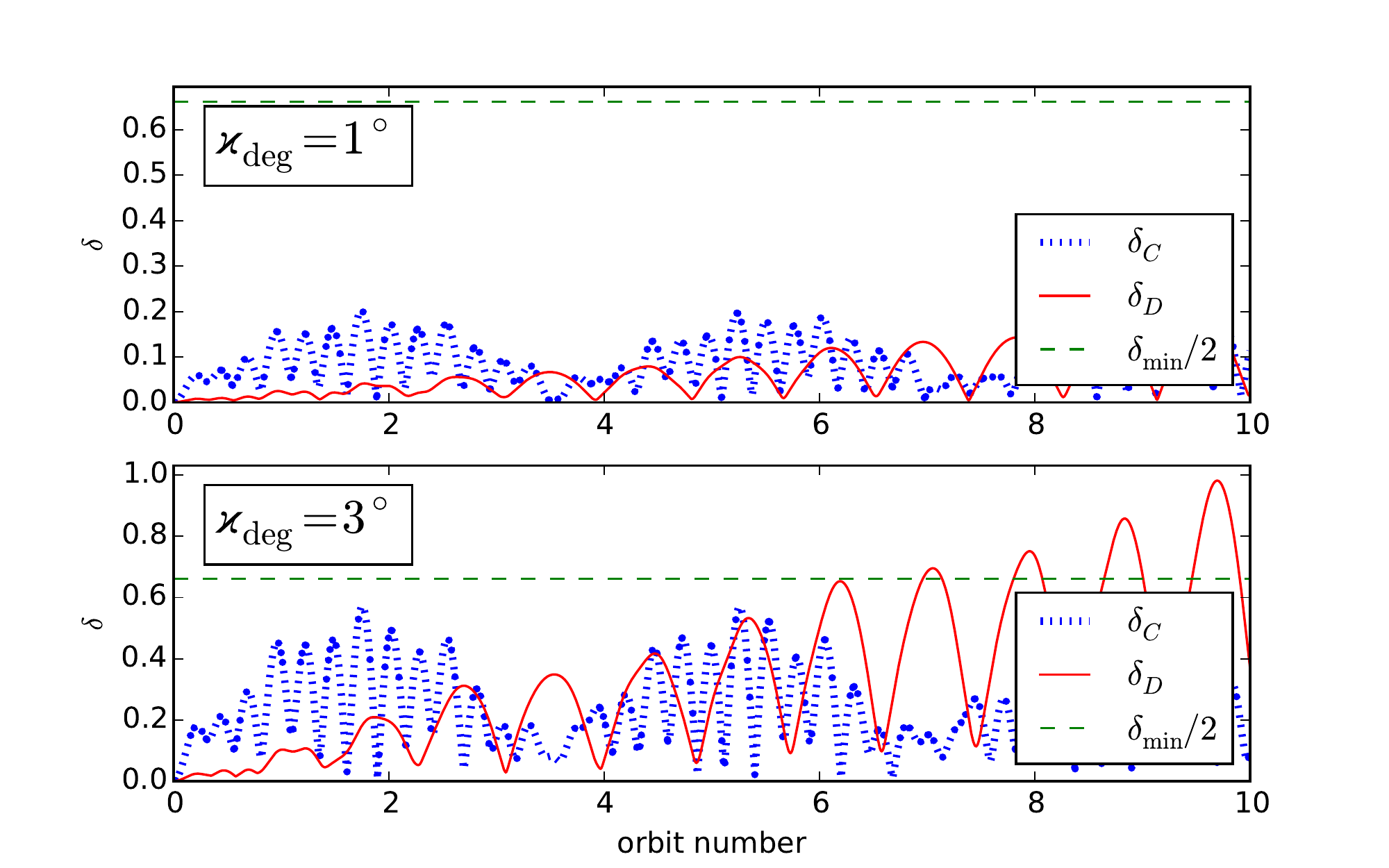}
    \caption{Deviations for Type II and $N=2$. In this case $\delta_{\rm min}=1.32$. At $\varkappa_{\rm deg}=1^\circ$ stability condition \eqref{eq:collision} holds with a large margin for the whole simulation interval of 10 orbital periods, while at $\varkappa_{\rm deg}=3^\circ$ it breaks down after six orbital periods.}
    \label{fig:devTypeII_N2}
\end{center}
\end{figure}

Results in the case of Type II formation at $N=2$ (Fig.~\ref{fig:devTypeII_N2}) are on the whole similar to those obtained in the first case (Type I, $N=3$). The deviations are in fact now higher than in that case, but this is somewhat compensated by the larger $\delta_{\rm min}$, so that the no-collision condition at $\varkappa_{\rm deg}=3^\circ$ again holds about up to the half of the simulation interval.

\subsection{Frequency adjustment by tuning the mass ratio}
The specific pattern of deviation growth observed in the above examples (linear growth in time with superimposed periodicity) strongly suggests that this growth is largely due to the gradual shift of the trajectories occurring because the ratio of the true $x$ and $y$ system's frequencies does not exactly match the approximate value $\alpha_0=N m_D / m_C = 8$. We can expect to negate these shifts by adjusting the mass ratio:
\[
\frac{N m_D}{m_C}=\alpha_0-\Delta\alpha.
\]

We find the appropriate $\Delta\alpha$ numerically, by minimizing the maximum deviation of deputy satellites for 10 orbital periods:
\[
\max_{t\in[0,10T_0]}\delta_D(t) \rightarrow \min_{\Delta\alpha}
\]

\begin{figure}[htb]
\begin{center}
\includegraphics[width=0.8\textwidth, keepaspectratio, clip, trim=0mm 7mm 0mm 0mm]{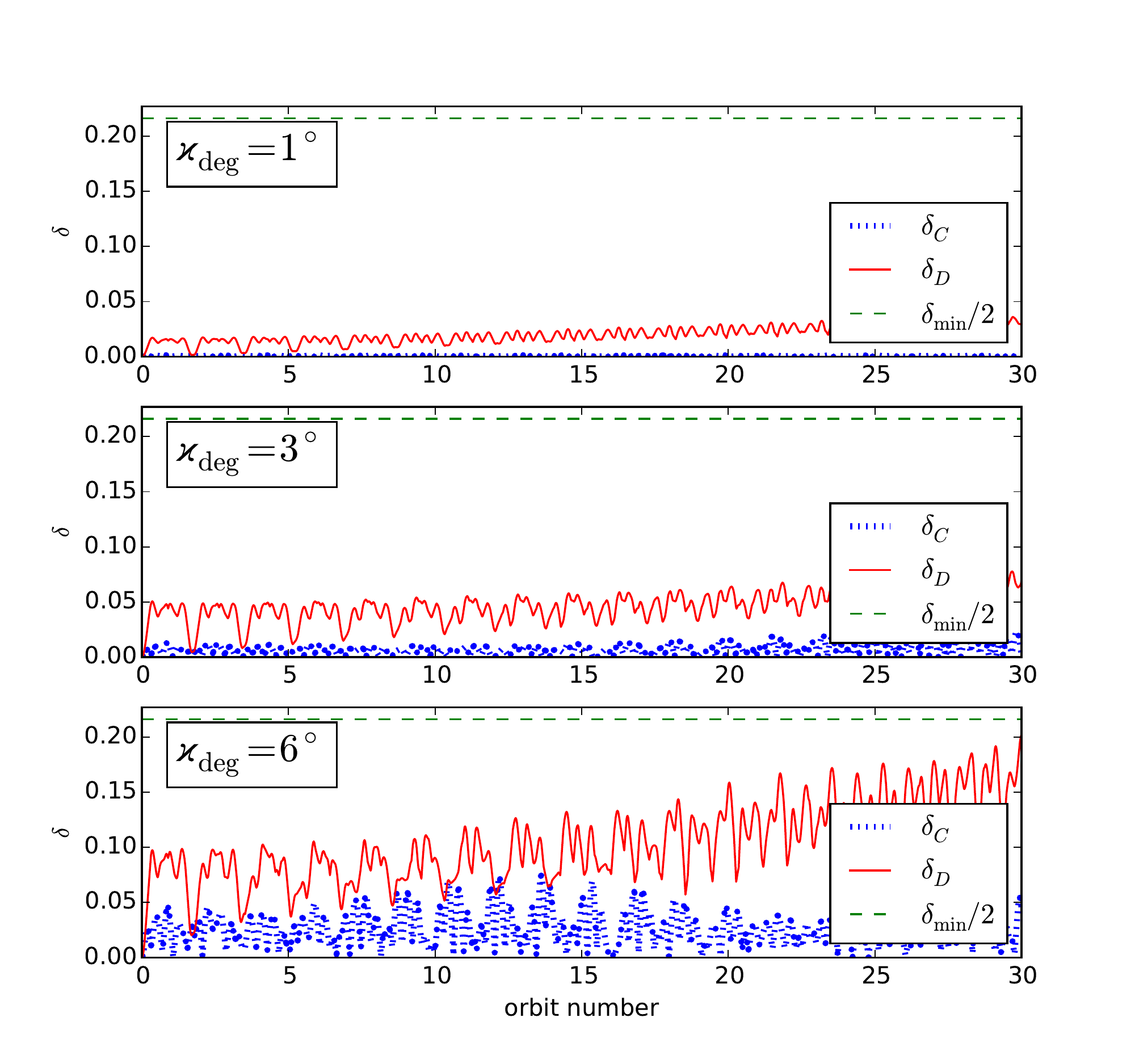}
\end{center}
    \caption{Deviations in the Type I, $N=5$ formation with the optimally adjusted satellite mass ratio. Stability condition \eqref{eq:collision} holds for the whole simulation interval of 30 orbital periods for all three considered values of $\varkappa_{\rm deg}$.}
    \label{fig:relDeviations_opt}
\end{figure}

We perform these optimizations at different angles $\varkappa_{\rm deg}$ for the second, most stable configuration from the previous section (Type I, $N=5$). The results are shown in Table 2. The results clearly show a big improvement over the earlier results obtained without mass ratio adjustment. In Fig.~\ref{fig:relDeviations_opt} the deviations of the adjusted systems with different values of $\varkappa_{\rm deg}$ are plotted for 30 orbital periods.

\begin{table}
\begin{center}
\begin{tabular}{ccccc}
\toprule
$\varkappa_{\rm deg}$ & $\Delta\alpha_{\rm opt}$ & $\max\limits_{t\in[0,10T_0]}\delta_D(t)$ & $\max\limits_{t\in[0,10T_0]}\delta_D(t)$ & $\max\limits_{t\in[0,30T_0]}\delta_D(t)$\\
 &  & w/o adjustment & with adjustment & with adjustment\\
\midrule
$1^\circ$ & 0.021 & 0.0408 & 0.0210& 0.0375  \\
$2^\circ$ & 0.039 & 0.0766 & 0.0368& 0.0530 \\
$3^\circ$ & 0.071 & 0.1230 & 0.0541& 0.0814  \\
$4^\circ$ & 0.116 & 0.1780 & 0.0702& 0.1280  \\
$5^\circ$ & 0.169 & 0.2600 & 0.0874& 0.1650   \\
$6^\circ$ & 0.230 & 0.3700 & 0.1110& 0.2020  \\
\bottomrule
\end{tabular}
\caption{Results of the mass ratio adjustment experiments. At each angle $\varkappa_{\rm deg}$ we numerically find the optimal adjustment $\Delta\alpha_{\rm opt}$. The maximum deviations of the main satellite from the vertical without adjustment or with the optimized adjustment are shown in the third and fourth columns, respectively. The last column shows results of simulations with adjustment spanning 30 orbital periods.	}
\end{center}
	
\end{table}

\section{Conclusions}
We found that subject to appropriate choice of ``hub-and-spoke'' system parameters deputy satellites can move along Lissajous curves so that the system stays in free motion (i.e. no fuel is consumed in the nominal operation mode). The main satellite is in a state of relative equilibrium on the local vertical, passing through the system CoM.

Our analysis shows the existence of rather nontrivial relations between the system's parameters that, when satisfied, produce a well-balanced system without collisions between the deputy satellites or tethers. Certain configurations determined by the proper choice of parameters allow placing an additional satellite at the center of the ``hub-and-spoke'' system without collisions with the other parts of the system.

Of course, the intricate way in which the deputy satellites change their positions is challenging for technical implementation. In particular, the tethers must be attached to the main satellite so as not to intertwine. Analysis of the tether entanglement in terms of homotopy provides another set of constraints that the system parameters should satisfy.

One of the most curious results in this paper is the study of nonlinear effects, which allows to formulate yet another set of conditions for the systems parameters to cancel out the second order corrections in the systems' equations of motion. Our numerical experiments with two parameter sets (for 5 deputy satellites) satisfying these conditions corroborate the theoretically predicted system's stability

The application of the proposed motion pattern is limited to small angular deviations of the tethers from the local vertical. The possibility to extend it to large deviations from the vertical could be a subject of further investigation.

\section*{Acknowledgment}
The authors first conceived the idea of the motion described in this paper during the dynamical analysis of the rotating multi-tethered satellite system \cite{A2015}, and the authors sincerely thank their collaborators Didier Alary, Kirill Andreev, Pavel Boyko, Elena Ivanova, and Cyrille Tourneur for the warm and stimulating atmosphere of that study. The work of one of the authors (DY) on the present paper was supported by Russian Science Foundation (project 14-50-00150).

\begin{appendices}
\section{Proof of Proposition 1}\label{sec:A}

\noindent Statements a) and d) are very simple, so we only provide proofs for b) and c).

b) Suppose that the satellite $i$ collides with the satellite $j$ at the time moment $\tau$. The condition $x_i(\tau)=x_j(\tau)$ admits two series of solutions:

\begin{enumerate}
\item  $2\pi p\left(\tau +\frac{i}{N}\right)+\varphi_x=
2\pi\left[p\left(\tau+\frac{j}{N}\right)+n_x\right]+\varphi_x$
with some integer $n_x$, that is
\[\frac{i-j}{N}=\frac{n_x}{p}.\]

\item $2\pi p\left(\tau +\frac{i}{N}\right)+\varphi_x=
2\pi\left[\frac{1}{2}+n_x-p\left(\tau+\frac{j}{N}\right)\right]-\varphi_x$
with some integer $n_x$, that is
$$
\tau =\frac{1}{2}\left[\frac{1}{p}\left(\frac{1}{2}+n_x-\frac{\varphi_x}{\pi}\right)-
\frac{i+j}{N}\right].
$$
\end{enumerate}

The condition $y_i(\tau )=y_j(\tau )$  has similar series but with $q$ instead of $p$, $\varphi_y$ instead of $\varphi_x$ and $n_y$ instead of $n_x$.

Thus there are four possibilities for the satellites to collide: an element in series 1 or 2 for the coordinate $x$ must occur simultaneously with an element in series 1 or 2 for the coordinate $y$. Let us deal with these cases one by one.

\begin{enumerate}
\item  Series 1 for $x$ and series 1 for $y$ -- impossible because $p$ and $q$ are co-prime and $-N<i-j<N$.

\item  Series 2 for $x$ and series 1 for $y$. Series 1 for $y$ has non-trivial solutions iff $N$ and $q$ are not co-prime. If a non-trivial solution exists, the time moment $\tau$ is found from series 2 for $x$. Thus in this series collisions occur iff $N$ and $q$ are not co-prime.

\item  Series 1 for $x$ and series 2 for $y$. Similarly, collisions occur iff $N$ and $p\ $ are not co-prime.

\item  Series 2 for $x$ and series 2 for $y$. Equating $\tau$ from both series, we obtain
$$
\frac{1}{p}\left(\frac{1}{2}+n_x-\frac{{\varphi }_x}{\pi }\right)=
\frac{1}{q}\left(\frac{1}{2}+n_y-\frac{{\varphi }_y}{\pi }\right)
$$
whence
\[n_xq{\rm -}n_yp{\rm =}\frac{{q\varphi }_x-{p\varphi }_y}{\pi }+\frac{p-q}{2}.\]
\end{enumerate}

Since $p$ and $q$ are co-prime, with $n_x,n_y$ running over all possible integers, the left side also runs over all possible integers, i.e. collisions in this case happen iff the right side is an integer.

c) Like in the proof of b), we obtain two series of relations from the condition $x_i{\rm (}\tau {\rm )=}x_{{\rm j}}{\rm (}\tau {\rm )}$. However, the first series takes the form
\[
\frac{i-j}{N}{\rm =}n_x
\]
and has no non-trivial solutions, because ${\rm -}N{\rm <}i-j{\rm <}N$. This leaves the second series, which has the form
\[
\tau =\frac{1}{2p}\left(\frac{1}{2}+n_x-\frac{\varphi_x}{\pi}-
\frac{i+j}{N}\right).
\]

Like before, we equate $\tau $ from the series for $x$ and $y$ and obtain
\[
\frac{1}{p}\left(\frac{1}{2}+n_x-\frac{\varphi_x}{\pi}-
\frac{i+j}{N}\right)=
\frac{1}{q}\left(\frac{1}{2}+n_y-\frac{\varphi_y}{\pi}-
\frac{i+j}{N}\right)
\]
whence
\[N\left(qn_x-pn_y\right)-\left(q-p\right)\left(i+j\right)=\left(\frac{{q\varphi }_x-{p\varphi }_y}{\pi }+\frac{p-q}{2}\right)N.\]

As $n_x{\rm ,\ }n_y$ run over all integer values, the expression $qn_x{\rm -}pn_y$ also takes all integer values. If $N\ge 3$, the sum $i+j$ for all possible pairs of different numbers from 1 to $N$ assumes all possible integer values modulo $N$. Thus the expression in the left side takes all possible integer values divisible by the greatest common divisor of $N$ and $q-p$. On the other hand, if $N=2$, then $i+j=3$ and the above relation simplifies to the requirement that
\[
\frac{{q\varphi }_x-{p\varphi }_y}{\pi }
\]
be integer.

\section{Proof of Proposition 2}\label{sec:B}

a) For a closed curve ${{\rm (}x\left(\tau \right),y(\tau ))}_{\tau \in [0,1]}$ not containing the origin, the number $w$ of its rotations about the origin can be written as
\[w=\frac{1}{2}\sum_{\tau :y\left(\tau \right)=0}{{\rm sign}\left(x\left(\tau \right)\dot{y}\left(\tau \right)\right)},\]
assuming that the intersections of the curve with the $x$ coordinate axis are non-degenerate. We will apply this formula to
\[
x_{ij}\left(\tau \right)=x_j\left(\tau \right)-x_i\left(\tau \right),\ y_{ij}\left(\tau \right)=y_j\left(\tau \right)-y_i\left(\tau \right).
\]
We consider separately the two types of formations.

\textit{Type I.} In this case
\[x_{ij}\left(\tau \right)=2x_0{\sin  \frac{\pi p(j-i)}{N}}{\cos  \left[2\pi p\left(\tau +\frac{j+i}{{\rm 2}N}\right)+{\varphi }_x\right],}\]
\[y_{ij}\left(\tau \right)=2y_0{\sin  \frac{\pi q(j-i)}{N}}{\cos  \left[2\pi q\left(\tau +\frac{j+i}{{\rm 2}N}\right)+{\varphi }_y\right],\ }\]
so that
\[w_{i,j}={\rm sign}\left(x_0y_0\right){\rm sign}\left[\sin \frac{\pi p(j-i)}{N}\right]
                                      {\rm sign}\left[\sin \frac{\pi q(j-i)}{N}\right]w_*=\pm w_*,\]
where $w_*$ is the winding number of the Lissajous curve
\[x_*(\tau )={\cos  \left(2\pi p\tau +{\varphi }_x\right),\ }\ y_*(\tau )={\cos  \left(2\pi q\tau +{\varphi }_y\right)\ }.\]
We have
\[
y_*\left(\tau \right)=0\ {\rm for}\
\tau =\frac{1}{2q}\left(s-\frac{{\varphi }_y}{\pi} - \frac{1}{2}\right),\ s=1,..2q.
\]

Moreover, with this choice of $\tau $ we have $\dot{y}_*\left(\tau \right)>0$ for even $s$ and $\dot{y}_*\left(\tau \right)<0$ for odd $s$. Applying the formula for the winding number, we obtain
\[
w_*=\frac{1}{2}
\sum^{2q}_{s=1} (-1)^s
{{\rm sign}\left[{\cos  \left(\frac{\pi ps}{q} + \varphi_*\right)\ }\right]}
\]
with some (unimportant) phase constant $\varphi_*$.

Now consider separately the cases when both $p$, $q$ are odd and when one of them is even.

Let both $p$, $q$ be odd. Since $p$, $q$ are co-prime, the values $(\pi ps/q){\rm mod}\ 2\pi $ run over the values in the set $Z_{2q}=\left\{\pi r/q \right\}^{2q-1}_{r=0}$ as $s$ runs over $1,2,..,2q$. Then the quantity
${\cos  \left( \frac{\pi ps}{q}+{\varphi }_*\right)\ }$
takes equally many positive and negative values as $s$ runs over $1,2,..,2q$. Now, if $s$ runs only over even values $2,..,2q$, then $(\pi ps/q){\rm mod}\ 2\pi $ runs over the values in the set $Z_q=\left\{2\pi r/q\right\}^{q-1}_{r=0}$. Then, since $q$ is odd, the number of positive and negative values taken by ${\cos  \left(\frac{\pi ps}{q}+{\varphi }_*\right)\ }$, as $s$ runs only even values $2,..,2q$, differs by 1. It follows that
\[
w_*=\frac{1}{2}\left(\pm 1-\left(\mp 1\right)\right)=\pm 1.
\]

Now let one of $p$, $q$ be even; without loss of generality we can assume that it is $q$. We can then repeat the above argument, but since $q$ is now even, we conclude that ${\cos  \left(\frac{\pi ps}{q}+{\varphi }_*\right)\ }$ takes equally many negative and positive values as $s$ runs over the even or odd subset of $1,2..,2q$. It follows that $w_*=\frac{1}{2}\left(0-0\right)=0$.

We have thus proved statement a) for Type I formation.

\textit{Type II.} In this case
\[
x_{ij}\left(\tau \right)=2x_0\sin \frac{\pi(j-i)}{N}\cos  \left[2\pi\left( p\tau +\frac{j+i}{2N}\right)+{\varphi }_x\right],
\]
\[
y_{ij}\left(\tau \right)=2y_0\sin \frac{\pi(j-i)}{N}\cos  \left[2\pi\left( q\tau +\frac{j+i}{2N}\right)+{\varphi }_y\right],
\]
which leads to
\[
w_{i,j}=\frac{{\rm sign}\left(x_0y_0\right)}{2}
\sum^{2q}_{s=1} (-1)^s
{\rm sign}\left(
\cos  \left[\pi\left(\frac{ps}{q}+\frac{i+j}{N}\frac{q-p}{q}\right)+{\varphi }_*\right]
\right)
\]
with some phase constant ${\varphi }_*$ not depending on $i$, $j$. The statement a) then follows just like in case of Type I.

b) From the proof of a), $w_{i,j}=w_*{\rm sign}\left(x_0y_0\right){\rm sign}(f_{j-i})$, where we set
\[
f_s=\sin  \left(\frac{\pi ps}{N}\right)\sin  \left(\frac{\pi qs}{N}\right).
\]
Note first that if $q-p$ is divisible by $2N$, then $f_s=\sin^2 \left(\pi p s/N\right)>0$ for all $s$ not divisible by $N$ (recall that $N$ and $p$ are co-prime by the assumption of no collisions). Accordingly, $w_{i,j}=w_*{\rm sign}\left(x_0y_0\right)$ for all $i,j=1,..,N$ with $i\ne j$. Similarly, $w_{i,j}=-w_*{\rm sign}\left(x_0y_0\right)$ if $q+p$ is divisible by $2N$.

Now suppose that neither of $q-p,\ q+p$ is divisible by $2N$. We need to show that $f_s$ takes both positive and negative values as $s$ runs over $1,2,..,N-1$. Since $f_N=f_0=0$ and $f_{2N-s}=f_s$, it suffices to show that $\sum^{2N-1}_{s=0}{f_s}=0.$ But that immediately follows from the hypothesis and the identity
\[
f_s=\frac{1}{2}\left[\cos  \frac{\pi s(p-q)}{N}-
                     \cos  \frac{\pi s(p+q)}{N}\right].
\]
The presence of strong entanglement follows from the presence of different winding numbers, since under a weak homotopy the winding number (possibly nonzero) must be the same for all pairs.
\end{appendices}

\end{document}